\documentclass[12pt]{article}
\usepackage{amsmath,amsthm,amssymb}

\newtheorem{theorem}{Theorem}

\newtheorem{lemma}{Lemma}
\newtheorem{definition}{Definition}

\newcommand{\ga}{\alpha}
\newcommand{\gb}{\beta}
\newcommand{\gc}{\gamma}
\newcommand{\gd}{\delta}
\newcommand{\gl}{\lambda}
\newcommand{\gs}{\sigma}
\newcommand{\hess}{\nabla^2}
\newcommand{\gra}{\nabla}
\newcommand{\de}{\partial}
\newcommand{\bpf}{\begin{proof}}
\newcommand{\epf}{\end{proof}}
\newcommand{\beq}{\begin{equation}}
\newcommand{\eeq}{\end{equation}}

\begin{document}
\title{Boundary Value Problems for some Fully Nonlinear
  Elliptic Equations}
 \vskip 1em
 \author{Szu-yu Sophie Chen}
 \date{December 5, 2005}
\maketitle


Let $(M, g)$ be a compact Riemannian manifold of dimension $n \geq
3$ with boundary $\de M$. We denote the Ricci curvature, scalar
curvature, mean curvature, and the second fundamental form by Ric,
R , h, and $L_{\alpha \beta},$ respectively.

The Yamabe problem for manifolds with boundary is to find a
conformal metric $\hat{g} = e^{-2u} g$ such that the scalar
curvature is constant and the mean curvature is zero. The boundary
is called umbilic if the second fundamental form $L_{\ga \gb}=
\mu_g g_{\ga \gb}$. For example, a totally geodesic boundary is
umbilic with zero principal curvatures. In \cite{Es92}, it was
proved by Escobar that for locally conformally flat compact
manifolds with umbilic boundary (and some other cases), the Yamabe
problem is solvable.

As for the nonlinear version of the Yamabe problem, we consider
the Schouten tensor defined as
$$ A_g =\frac{1}{n-2} ( Ric - \frac{R}{2(n-1)} g ) .$$
Note that $\mathrm{tr} A_g = \frac{1}{2(n-1)} R.$ The Schouten
tensor comes naturally from  curvature decomposition
 $$Riem = \mathcal{W} + A \odot g,$$ where the Weyl tensor $\mathcal{W}$
 is locally conformally invariant, and $\odot$ stands for the Kulkarni-Nomizu
 product. In dimension four, we have the
 following Chern-Gauss-Bonnet formula for closed manifolds:
$$ 32 \pi^2 \chi (M^4) = \int_{M^4} |\mathcal{W}|^2 + 16 \int_{M^4} \gs_2 (A_g),$$
where $\chi$ is the Euler characteristic and $\gs_2 (A_g)$ is the
second elementary symmetric function of the eigenvalues of $A_g.$
Since $\chi$ is a topological invariant and $\mathcal{W}$ is
locally conformally invariant, we have that $\int_M \gs_2 (A_g)$
is a conformal invariant. For closed four-manifolds,
Chang-Gursky-Yang \cite{CGY02a} proved that if the Yamabe constant
and $\int_M \gs_2 (A_g)$ are both positive, then we can find a
conformal metric $\hat g$ such that $\gs_2(A_{\hat g})$  is
constant. For locally conformally flat closed manifolds, Li-Li
\cite{LL03} proved that if $\gs_i(A_g) > 0, 1 \leq i \leq k$ for
some $k \geq 2,$ then we can find a conformal metric $\hat g$ such
that $\gs_k(A_{\hat g})$ is constant. See also Guan-Wang \cite{GW03a}
for an independent work of the above result. For closed manifolds which are not locally conformally
 flat,  Gursky-Viaclovsky \cite{GV05} proved that if $\gs_i(A_g) > 0$,
 for $1 \leq i \leq k$ and $2k > n,$ then we can
 find a conformal metric $\hat g$ such that $\gs_k(A_{\hat g})$ is
 constant.


In this paper, we study the nonlinear version of Yamabe problem
for manifolds with boundary. Before introducing the problem, we
need the following definitions:
\begin{definition} Let $W$ be a matrix with eigenvalues $\gl_1,\cdots,\gl_n.$
 \; \, Then  $\gs_k (W) =$\; $ \sum_{i_1<\cdots< i_k} \gl_{i_1}\gl_{i_2}\cdots\gl_{i_k}$
 for $k \leq n$ is called the
kth elementary symmetric function of the eigenvalues of $W$.
Denote $\gs_0 = 1.$ For example, $\gs_1 = \gl_1 + \cdots + \gl_n =
\mathrm{tr} \, W$ and $\gs_n = \gl_1 \cdots \gl_n = \det W.$
\end{definition}
The elementary symmetric functions are special cases of hyperbolic
 polynomials introduced by Garding \cite{Gar59}, which have nice
 properties in associated cones.
 \begin{definition}
 The set $\Gamma^+_k = \{$ the connected component of $\gs_k(\gl)>
0$ which contains the identity $\}$ is called the positive
$k$-cone.
 Equivalently, it is showed in \cite{Gar59} that
 $\Gamma^+_k = \{ \gl \,: \gs _i (\gl)> 0 , 1 \leq i \leq k \}$
 is an open convex cone with vertex at the origin, e.g.,
 $\Gamma^+_1 = \{\gl : \gl_1 + \cdots \gl_n > 0\}$ and $\Gamma^+_n =
 \{\gl : \gl_i > 0 , 1 \leq i \leq n \}.$ The following is the nested
 relation
  $$ \Gamma^+_1 \supset \Gamma^+_2 \supset \cdots \supset \Gamma^+_n.$$
 Denote $W \in \Gamma^+_k$ if the eigenvalues $\gl(W) \in \Gamma^+_k.$
\end{definition}
Suppose that the boundary is umbilic. Our goal is to find a
conformal metric $\hat g = e^{-2u}g$ such that  $\gs_k(A_{\hat
g})$ is constant and the boundary is totally geodesic.  We now
describe a class of locally conformally flat compact manifolds of
dimension $n\geq 3$ with boundary, for which we give an
affirmative answer to the question. Under the conformal change of
the metric $\hat g = e^{-2u}g$, we denote the curvature tensors in
the new metric by a \emph{hat} (For example, $\hat A, \hat L$ and
$ \hat \mu$). The Schouten tensor $\hat A$ satisfies
 \beq \label{e:schten}
 \hat A= \hess u + du \otimes du - \frac{1}{2}|\gra u|^2 g +
 A_g,
 \eeq
  where the derivatives are covariant derivatives with respect to
 the background metric $g.$ The second fundamental form satisfies
 $$\hat L e^{u}= \frac{\de u}{\de n} g + L_g,$$
 where $n$ is the
 unit inner normal with respect to $g$ on the boundary. Note that umbilicity is
 conformally invariant. Thus, it is natural to consider the class of
 manifolds with umbilic boundary.  When the boundary is umbilic,
 the above formula becomes
  $$ \hat \mu e^{-u} = \frac{\de u}{\de n} + \mu_g.$$

 If we view $\hat A$ as a $(0,2)$-tensor
 in the new metric $\hat g,$ then  $\gs_k(\hat A): \equiv \gs_k(\hat g^{-1}
 \hat A),$ where $\hat g^{-1}$ is the induced inverse tensor of
 the metric tensor $\hat g.$ On the other hand, by formula (\ref{e:schten})
 we can also view $\hat A$ as a $(0,2)$-tensor in the background
 metric $g.$ Using this notation, the problem becomes to consider
 the following equation:
 \begin {equation} \label {e:sigma}\left\{  \begin{array}{ll}
 \sigma_k^{\frac{1}{k}} (\hess u + du \otimes du - \frac{1}{2}|\gra u|^2 g +
 A_g)= e^{-2u} & in \, M  \\
 \frac{\de u}{\de n} + \mu_g = 0& on \,\partial M.
  \end{array}\right .
  \end {equation}


\begin{theorem} \label{t:sigma}
Suppose $(M, g)$ is a locally conformally flat compact manifold of
dimension $n\geq 3$ with umbilic boundary. If $A_g \in \Gamma^+_k$
for $k \geq 2$, then there exists a smooth solution $u$ of
(\ref{e:sigma}). In other words, there is a conformal metric $\hat
g = e^{-2u}g$ such that $\gs_k(\hat A)= 1$ and the boundary is
totally geodesic.
\end{theorem}

We will prove a more general result than Theorem~\ref{t:sigma}.
Consider the equation
\begin {equation} \label {e:main}\left\{  \begin{array}{ll}
 F(\hess u + du \otimes du - \frac{1}{2}|\gra u|^2 g +
 A_g)=  e^{-2u} & in \, M  \\
 \frac{\de u}{\de n} + \mu_g = 0& on \,\de M,
  \end{array}\right .
  \end {equation}
 where $F$ satisfies some structure conditions
 listed below. Equation (\ref{e:main}) means that
 we apply $F$ to the eigenvalues of the matrix (or $(1,1)$-tensor)
 $g^{-1}(\hess u + du \otimes du - \frac{1}{2}|\gra u|^2 g +
 A_g ).$ Now we give structure conditions for $F.$
 Let $\Gamma$ be an open convex cone in $\mathbb{R}^n$ with vertex at the origin
 satisfying $\Gamma^+_n \subset \Gamma \subset \Gamma^+_1.$
 Suppose that $F(\gl)= F(\gs_1(\gl),\cdots,\gs_n(\gl)) \in C^{\infty}(\Gamma) \cap C^0(\overline{\Gamma})$
 is a homogeneous symmetric
 function of degree one normalized with $F(e) = F(1,\cdots,1) = 1.$
 Assume that $F = 0$ on $\de \Gamma$ and $F$ satisfies the following in $\Gamma:$

 (S0) $F$ is positive;\par
 (S1) $F$ is concave (i.e., $\frac{\de^2 F}{\de \gl_i \de
 \gl_j}$ is negative semi-definite);\par
 (S2) $F$ is monotone (i.e., $\frac{\de F}{\de \gl_i}$ is
 positive);\par
 (S3) $\frac{\de F}{\de \gl_i}\geq \epsilon \frac{F}{ \gs_1},$ for some constant
     $\epsilon > 0,$ for all $i.$ \par
 In some case, we need an additional condition: \par
  (A) $\sum_{j \neq i} \frac{\de F}{\de \gl_j} \leq \rho \frac{\de F}{\de \gl_i},$
      for some $\rho > 0,$ for all $\gl \in \Gamma$ with $\gl_i \leq 0.$
 \vskip 1em
 An easy example is $F = \frac{1}{n}(\gl_1 + \cdots + \gl_n)$ with $\Gamma =
 \{\gl : \gl_1 + \cdots + \gl_n > 0 \}.$ Condition (S1) is used in most
 elliptic theories. Condition (S2) is the actual ellipticity. It
 is an elementary fact that if $F$ is a symmetric function of
 eigenvalues, then $\frac{\de F}{\de \gl_i}> 0 $ for all $i$ if and only if
 $F^{ij} :\equiv \frac{\de F}{\de W_{ij}}$ is positive definite.
 Condition (S3) was before in \cite{CGY02a}.

 \begin{theorem}\label{t:main}
 Suppose $(M, g)$ is a locally conformally flat
compact manifold of dimension $n\geq 3$ with umbilic boundary. Let
$F$ satisfy the structure conditions (S0)-(S3) in a corresponding
cone $\Gamma.$ If $A_g \in \Gamma,$ then there exists a smooth
solution $u$ of (\ref{e:main}).
 \end{theorem}

 In Section~\ref{s:sigma} below, we will show that $\binom{n}{k}^{-\frac{1}{k}}
  \gs^{\frac{1}{k}}_k$ satisfies the structure conditions
  (S0)-(S3)
  with $\epsilon = \frac{1}{k}$ in $\Gamma^+_k.$
 Hence, Theorem~\ref{t:main} implies Theorem~\ref{t:sigma}.

The next result concerns boundary estimates for equations more
general than (\ref{e:main}). Before stating the theorem, we
introduce some notations. In this paper, we use Fermi (geodesic)
coordinates in a boundary neighborhood, which means that we take
the geodesics in the normal direction parameterized  by arc length
from a local chart $(x_1, \cdots$ $, x_{n-1})$ on the boundary.
The metric is then expressed as $g = dx^n dx^n + g_{\ga \gb}
dx^{\ga} dx^{\gb} $.
 The Greek letters $\ga,\gb,\gc$ stand for the tangential direction
indices, $1\leq \ga,\gb,\gc < n,$ while the letters $i,j,k$ stand
for the full indices, $1 \leq i,j,k \leq n$. Define the half ball
in Fermi coordinates by $\overline{B}_r^+ =
  \{ x_n \geq 0, \sum_i x_i^2 \leq r^2\}$ and the segment on the
  boundary by $\Sigma_r =  \{ x_n = 0, \sum_i x_i^2 \leq r^2\}.$
All derivatives are covariant derivatives with respect to the
background metric $g$ unless otherwise noted.

The following boundary estimates are used in the proof of
Theorem~\ref{t:main}.

\begin{theorem} \label{t:bdy}
 Let $F$ satisfy (S0)-(S3) in a corresponding cone $\Gamma$
 and $g$ be a flat metric.
 Suppose that $\Sigma_r$ is umbilic with principal curvatures $\mu$ and $n$ is the unit
 inner normal with respect to $g.$
 Let $u \in C^4$ be a solution to the equation
 \begin {equation} \left\{  \begin{array}{ll}
 F(\hess u + du \otimes du - \frac{1}{2}|\gra u|^2 g
 )= f e^{-2u} & in \, \overline {B}^+_r  \\
 \frac{\de u}{\de n} + \mu =  \hat \mu\, e^{-u} & on \, \Sigma_r.
  \end{array}\right .
  \end {equation}
 Case(a). If $\hat \mu = 0,$ then
 $$\sup_{x \in \overline{B}^+_{\frac{r}{2}}} \,( |\gra u|^2 + |\hess u| ) \leq
   C ( 1 + \sup_{x \in \overline{B}^+_r} e^{ -2 u}),$$
  where $C$ depends on $ r, n, \epsilon, \mu, \| f\|_{C^2(\overline {B}^+_r)}$ and
  $\inf_{\overline {B}^+_r} f.$\\
 Case(b).  Suppose that $F$ satisfies
 the additional condition (A) and $\Gamma_2^+ \subset \Gamma.$
  If $\hat \mu$ is a positive constant, then
  $$\sup_{x \in \overline{B}^+_{\frac{r}{2}}} \,( |\gra u|^2 + |\hess u| ) \leq
   C, $$
  where $C$ depends on $ r, n, \epsilon, \rho, \mu, \hat \mu, \inf_{\overline {B}^+_r} u,
   \| f\|_{C^2(\overline {B}^+_r)}$ and $\inf_{\overline {B}^+_r} f.$
 \end{theorem}
   In Section~\ref{s:sigma} below, we  further show that $ \binom{n}{k}^{-\frac{1}{k}}
  \gs^{\frac{1}{k}}_k$ satisfies the additional condition (A)
  with $\rho = (n- k).$ Thus, $ \binom{n}{k}^{-\frac{1}{k}}
  \gs^{\frac{1}{k}}_k$ for $k \geq 2$ is an example of case (b).

 The Dirichlet problems for fully nonlinear elliptic equations have been
 extensively studied, for example, by Caffarelli-Nirenberg-Spruck \cite{CNS-1}, \cite{CNS-3}
 and by Trudinger \cite{Tr90}.
 Such problems  for the Schouten tensor equations are studied by
 Guan \cite{Gb05}. On the other hand, the Neumann problems for fully
 nonlinear elliptic equations are not yet well studied. The problem we
 proposed here comes from natural geometrical setting. It would be
 an interesting problem whether we can consider other
 Monge-Ampere-type equations.

 The idea of proof of Theorem~\ref{t:main} is to deform the Yamabe metric
for  manifolds with boundary to the one satisfying the equation
(\ref{e:main}). The similar idea has already appeared in
\cite{LL03} and \cite{GV03} for closed manifolds. We will show
that, to avoid the bubbling phenomenon, if a manifold is not
conformally equivalent to  hemispheres, we have a priori
estimates.  Hence by degree theory argument we obtain a solution.
The proof of boundary $C^0$ estimates follows closely that of
Li-Li \cite{LL03}, while we still need to prove a revised version
of the work by Schoen-Yau \cite{SY88}, which turns out to be a
crucial element. As for $C^2$
 estimates, local $C^2$ estimates are previously
proved by Chang-Gursky-Yang \cite{CGY02}, Guan-Wang \cite{GW03}
and Li-Li \cite{LL03} in different cases. Recently, a simplified
proof of local $C^2$ estimates is derived by Chen \cite{Chen05}
and applied to a large class of equations. To prove
Theorem~\ref{t:bdy}, we will use an idea in that work to derive
boundary $C^2$ estimates directly from boundary $C^0$ estimates,
which is the main part of this paper.

The above results extend to manifolds with boundary which are not
locally conformally flat. In a subsequent paper \cite{Chen05b}, we
study boundary value problems associated to some integral
invariants on manifolds with boundary.

The paper is organized as follows. We start with some background
in Section~\ref{s:sigma}.  In Sections~\ref{s:main} and
\ref{s:bdy}, we give the proofs of Theorem~\ref{t:main} and
\ref{t:bdy}, respectively.

 \vskip 1em
 \textbf{Acknowledgment:} The author would like to thank Alice
 Chang for her constant support during the author's graduate
 education at Princeton University.


\section{Background}\label{s:sigma}
 We give some basic facts about homogeneous symmetric functions.
\begin{lemma}(see \cite{Chen05}).\label{l:sym}
 Let $\Gamma$ be an open convex cone with vertex at the origin
 satisfying $\Gamma^+_n \subset \Gamma$ ,and let
 $e = (1, \cdots, 1)$ be the identity.
 Suppose that $F$ is a homogeneous symmetric function of degree one
 normalized with $F(e)= 1,$ and that $F$ is concave in $\Gamma.$
 Then
 \begin{description}
 \item{(a)} $\sum_i \gl_i \frac{\de F(\gl)}{\de \gl_i} = F(\gl), \quad$ for $\gl \in
 \Gamma.$
 \item{(b)} $\sum_i \frac{\de F(\gl)}{\de \gl_i} \geq F(e) = 1, \quad$
 for $\gl \in \Gamma.$
 \end{description}
\end{lemma}
 Now we list further properties of elementary symmetric functions.
 \begin{lemma} \label{l:sigma}(see \cite{Gar59}, \cite{Mit70} and \cite{CNS-3}).
  Let $G =\gs_k ^{\frac{1}{k}}, k \leq n.$ Then\par
  \vskip 1em
  (a) $G$ is positive and concave in $\Gamma^+_k.$\par
  (b) $G$ is monotone in $\Gamma^+_k,$ i.e., the matrix
  $G^{ij} = \frac {\de G}{\de W_{ij}}$ is positive definite.\par
  (c) For $0 \leq l < k \leq n,$ the following is the Newton-MacLaurin inequality
    $$ k(n-l+1) \gs_{l-1} \gs_k \leq l(n-k+1) \gs_l \gs_{k-1}.$$
\end{lemma}
   Therefore, $S = \binom{n}{k}^{-\frac{1}{k}}\, G$
  satisfies the structure conditions (S0)-(S2) in $\Gamma^+_k.$
\vskip 1em
  We use the notation $\Lambda_i = (\gl_1, \cdots, \hat {\gl_i}, \cdots, \gl_n).$
 We will show that $S= \binom{n}{k}^{-\frac{1}{k}}\, G$ satisfies (S3) by using
 the following lemma:
 \begin{lemma}
  Let $n \geq 2.$ If $\gl \in \Gamma^+_k$ for some $1 \leq k \leq n,$ then
 $$\gs_{k-1}(\Lambda_i) \geq \frac{\gs_k (\gl)}{\gs_1 (\gl)} \qquad \forall i. $$
 \end{lemma}
 \bpf
 Since $\gl \in \Gamma^+_k,$ we have $\frac{\de \gs_l}{\de \gl_i} =
 \gs_{l-1} (\Lambda_i) > 0,$ for $1 \leq l \leq k,$ and thus
 $\Lambda_i \in \Gamma^+_{k-1} (\mathbb{R}^{n-1}).$
 On the other hand, by definition we have the identity $\gs_{k-1} (\Lambda_i) \gs_1(\gl)
 = \gs_{k-1} (\Lambda_i) \gl_i + \gs_{k-1} (\Lambda_i) \gs_1 (\Lambda_i).$

 \noindent Case (1): For $k=1,$ we get $\gs_{k-1}(\Lambda_i) = 1 = \frac{\gs_1 (\gl)}{\gs_1 (\gl)}.$

 \noindent Case (2): For $2 \leq k \leq n-1,$ by Lemma~\ref{l:sigma} (C),
   $(n-k) \gs_1 (\Lambda_i) \gs_{k-1} (\Lambda_i) \geq k (n-1) \gs_k (\Lambda_i).$
   If $\gs_k(\Lambda_i) \geq 0,$ then $$\gs_1 (\Lambda_i) \gs_{k-1} (\Lambda_i) \geq
    \frac{k (n-1)}{n-k} \gs_k (\Lambda_i) \geq \gs_k
    (\Lambda_i).$$
   If $\gs_k(\Lambda_i) < 0,$ then $$\gs_1 (\Lambda_i) \gs_{k-1} (\Lambda_i)
   > 0 > \gs_k (\Lambda_i).$$ Thus, in both cases, $\gs_{k-1} (\Lambda_i) \gs_1(\gl) \geq
    \gs_{k-1} (\Lambda_i) \gl_i + \gs_k (\Lambda_i) = \gs_k (\gl).$

 \noindent Case (3): For $k = n,$ we have
   $\gs_{n-1} (\Lambda_i) \gs_1 (\gl)  \geq \gs_{n-1} (\Lambda_i) \gl_i = \gs_n (\gl).$
 \epf
 As a consequence of the above lemma, $S = \binom{n}{k}^{-\frac{1}{k}}\, \gs_k ^{\frac{1}{k}}$
 satisfies (S3) with $\epsilon = \frac{1}{k}.$

 The next lemma shows that  $S = \binom{n}{k}^{-\frac{1}{k}}\, \gs_k
 ^{\frac{1}{k}}$ also satisfies the additional condition (A) with
 $\rho = (n-k).$
 \begin{lemma}
  For $1 \leq k \leq n-1 ,$ if $\gl \in \Gamma^+_k$ with $\gl_i \leq 0$ for some $i,$ then
  $$ \sum_{j \neq i} \frac{\de \gs_k (\gl)}{\de \gl_j} \leq (n- k) \frac{\de \gs_k (\gl)}{\de \gl_i}.$$
 \end{lemma}
 \bpf
   For $k=1,$ the above inequality is trivial since $\frac{\de \gs_1(\gl)}{\de \gl_j} =
   1$ for all $j.$
   For $k \geq 2$, we have
   \begin{eqnarray*}
   \sum_j \frac{\de \gs_k (\gl)}{\de \gl_j} &=& (n-k+1) \gs_{k-1}(\gl) =
  (n-k +1) (\gs_{k-1} (\Lambda_i)+ \gl_i  \gs_{k-2}(\Lambda_i))\\
     &\leq& (n-k+1) \gs_{k-1} (\Lambda_i) = (n-k+1) \frac{\de \gs_k
     (\gl)}{\de \gl_i}.
    \end{eqnarray*}
   By cancelling out $\frac{\de \gs_k (\gl)}{\de \gl_i}$ on both sides,
   the lemma is proved.
 \epf
 Suppose that $F$ satisfies (S0)-(S3) in $\Gamma.$
 It is useful to consider the following symmetric functions, which are
  introduced in \cite{LL03}.
\begin{definition} Let
 $F^t (\gl)= (t + n (1-t))^{-1} F(t \gl + (1-t) \gs_1 (\gl) e),$ for $0 \leq t \leq 1$
 in the cone $\Gamma^t = \{\gl: t \gl + (1-t) \gs_1(\gl) e \in \Gamma \}.$
\end{definition}
 We  show that $F^t$  satisfies (S0)-(S3) in $\Gamma^t.$ It is
 easy to see that $F^t$ is positive and concave. For monotonicity,
  $$(t + n (1-t)) \frac{\de F^t}{\de \gl_i} = t F_i + (1-t) \sum_j F_j \geq F_i > 0.$$
 As for (S3), $$\frac{\de F^t}{\de \gl_i} \geq \epsilon \frac{F^t(\gl)}{\gs_1(t \gl + (1-t) \gs_1(\gl) e)}
  (t + n(1-t)) = \epsilon \frac{F^t(\gl)}{\gs_1(\gl)}.$$
 Finally, if $F(\gl) = F(\gs_1(\gl), \cdots, \gs_n (\gl)),$ then
 $F^t (\gl) $ is a function of $\gs_1 (\gl), \cdots, \gs_n (\gl).$
 This is because $\gs_k(t \gl + (1-t) \gs_1(\gl) e),$ a homogeneous
 symmetric polynomial, is a function of $\gs_1(\gl), \cdots, \gs_n(\gl)$
 by elementary algebra.

The next lemma concerns  some important behaviors of solutions on
the boundary. As we mentioned in the introduction, in this paper
we use Fermi coordinates in a boundary neighborhood. Before
stating the lemma, we introduce a definition:
\begin{definition} (see \cite{Reilly}).
Let $P$ be a symmetric matrix. $T_k=\gs_k\,I - \gs_{k-1} P +
\cdots + (-1)^k P^k$ is called the $k$-th Newton tensor associated
with $P.$ We have that $\frac{\de \gs_k (P)}{\de P_{ij}} =
(T_{k-1})_{ij}.$
\end{definition}

\begin{lemma} \label{l:bdy} Let $F = F(\gs_1( g^{-1}\hat A), \cdots, \gs_n(g^{-1}\hat A)).$
 Suppose $g$ is flat and $\hat L_{\ga \gb} = \hat \mu \hat g_{\ga \gb}$ for some constant $\hat \mu$ near a
 boundary point $x_0.$ Then
 \begin{description}
 \item{(a)} $F^{\ga n} = 0$ at $x_0,$
 \item{(b)} $\hat A_{\ga \gb, n} = 2 \mu \hat A_{\ga \gb}
 - \hat \mu e^{-u} (\hat A_{\ga \gb} + \hat A_{nn} g_{\ga \gb})$ at $x_0.$
 \end{description}
\end{lemma}
 \bpf Since $g$ is flat, we have $\hat A = \hess u + du \otimes du - \frac{1}{2} |\gra u|^2 g.$
 We denote the covariant differentiation with respect to the new
 metric $\hat g$ by $\hat \gra.$ By the Codazzi equation
 $$
  \hat R_{\ga \gb \gc n} =  \hat \gra_{\gb }\hat L_{\ga \gc} -
 \hat \gra_{\ga} \hat L_{\gb \gc},
 $$ we have $\hat R_{\ga \gb \gc n} = 0$ because $\hat \mu$ is constant.
 Thus, we obtain $\hat R_{\ga n} = 0$ and $\hat A_{\ga n}= 0$ at $x_0.$
 To prove (a), since $F$ is a function of $\gs_i$, we only need to
 show that $\frac{\de \gs_i (g^{-1} \hat A)}{\de \hat A_{\ga n}} = (T_{i-1})_{\ga n} = 0$
for all $i.$ We prove it by induction. For $i= 1,$ by definition
$(T_1)_{\ga n} = \gs_1(g^{-1} \hat A) g_{\ga n} - \hat A_{\ga n},$
which equals to zero. For general $i,$ notice the recursive
relation $(T_i)_{\ga n} = \gs_i (g^{-1} \hat A) g_{\ga n} -
(T_{i-1})_{\ga j} \hat A_{j n}.$ Applying the induction hypothesis
gives $(T_i)_{\ga n} = - (T_{i-1})_{\ga \gb} \hat A_{\gb n} = 0.$

For (b), note that the boundary is umbilic. Thus, $u$ satisfies
 $\frac{\de u}{\de n} + \mu = \hat \mu e^{-u}$ on the boundary near $x_0.$
 Since $g$ is flat, by the Codazzi equation, $\mu$ is a constant.
 Notice that $\Gamma^n_{\ga n} = 0, \Gamma^n_{\ga \gb} = \mu g_{\ga \gb}$
 and $\Gamma^{\gb}_{\ga n} = - \mu \gd_{\ga \gb}.$
 Using the boundary condition, straightforward computations give us
 \beq \label{e:nga}
 u_{n \ga} = - \hat \mu e^{-u} u_{\ga} - \sum_j \Gamma^j_{\ga n} u_j = \mu u_{\ga} - \hat \mu u_{\ga}
 e^{-u},
 \eeq
 and
 \beq \label{e:ngagb} \begin{array}{ll}
 u_{\ga \gb n} &= (\mu - \hat \mu e^{-u}) (u_{\ga \gb} + \sum_j \Gamma^j_{\ga \gb} u_j)
  + \hat \mu u_{\ga} u_{\gb} e^{-u}
    - \sum_l \Gamma^l_{\gb n} u_{l \ga} - \sum_l \Gamma^l_{\ga \gb} u_{n l}\\
 &= (2 \mu -\hat \mu e^{-u}) u_{\ga \gb} - \mu u_{nn} g_{\ga \gb}
     + \hat \mu u_{\ga} u_{\gb} e^{-u} - \mu (-\mu + \hat \mu e^{-u})^2 g_{\ga \gb}.
 \end{array}
 \eeq
 Thus,
 \begin{eqnarray*}
 \hat A_{\ga \gb, n} &=& u_{\ga \gb n} + u_{\ga n} u_{\gb}
 + u_{\gb n} u_{\ga} - \sum_l u_l u_{l n} g_{\ga \gb}\\
  &=& 2 \mu ( u_{\ga \gb} +  u_{\ga} u_{\gb} - \frac{1}{2} |\gra u|^2
g_{\ga \gb}) - \hat \mu e^{-u} (u_{\ga \gb} + u_{\ga} u_{\gb} +
(-\sum_{\gc} u^2_{\gc} + u_{nn}) g_{\ga \gb}),
 \end{eqnarray*} which equals to $2 \mu \hat A_{\ga \gb}
 - \hat \mu e^{-u} (\hat A_{\ga \gb} + \hat A_{nn} g_{\ga \gb}).$
 \epf

 Remark: In above lemma, (b) can be proved in an another way. Since
 $g$ is flat, $\hat {\mathcal{W}}$ vanishes. Thus, by curvature
 decomposition $\hat R_{ijkl}$ can be written in terms of
 $\hat R_{ij}.$ Then using the Bianchi identity, we can compute $\hat A_{\ga \gb, n}.$

\section{Proof of Theorem~\ref{t:main}} \label{s:main}

 \bpf
 We deform the Yamabe metric to the one satisfying the equation
(\ref{e:main}). Define $F^t = (t + n (1-t))^{-1} F(t \gl + (1-t)
\gs_1 (\gl) e)$ in $\Gamma^t$ as in Section~\ref{s:sigma}. Let the
background metric $g$ be the Yamabe metric such that $R_g$ is a
positive constant and the boundary is totally geodesic. Thus, the
equation becomes the following: \beq \label {e:main'}\left\{
\begin{array}{ll}
 F^t(\hess u + du \otimes du - \frac{1}{2}|\gra u|^2 g +
 A_g)= e^{-2u} & in \, M  \\
 \frac{\de u}{\de n} = 0& on \,\partial M.
  \end{array}\right .
\eeq We will derive later a priori estimates for this path of
equations for $(M, g)$ not conformal equivalent to standard
hemispheres $(S^n_+, g_c),$ where $g_c$ is the standard metric on
spheres. The Leray-Schauder degree is defined similarly as in Li
\cite{Li89}. In our case, we just consider the space $\{u \in
C^{4, \ga} (M): \frac{\de u}{\de n} = 0$ on $\de M \}$  instead of
$\{u \in C^{4, \ga} (M)\}$ for most closed manifolds cases.
 Then by homotopy-invariance we obtain a solution at $t =1,$
 since at $t = 0$ the degree is nonzero.
 The fact that at $t=0$ the degree is nonzero is proved
 by Schoen \cite{Sch91} for the Yamabe problem on closed manifolds.
  In our case, $\frac{\de u}{\de n}= 0$ on $\de M$ so the boundary
 integral terms vanish in the computations in \cite{Sch91}. Thus,
 the result remains the same. The problem then reduces to
 establishing  a priori estimates.

 Suppose $F$ satisfies conditions (S0)-(S3). As
in the discussion in Section~\ref{s:sigma}, $F^t$ also satisfies
(S0)-(S3).  We drop $t$ without loss of generality in proving a
priori estimates. We denote the conformal equivalence relation by
$\cong.$

\noindent (1) $C^0$ estimates for $(M, g) \ncong (S^n_+, g_c).$

 Since the boundary is totally geodesic, it is natural to consider
 the doubling of the manifold $(M, g)$ and apply the $C^0$
 estimates on locally conformally flat
 closed manifolds. However, one problem is that we need to
 verify the doubling of the manifold still inherits a locally conformally flat
 smooth structure. Another problem is that the work by Schoen-Yau \cite{SY88}
 is for locally conformally flat smooth manifolds, which is a crucial element
 in the proof of $C^0$ estimates.  Thus, we need a revised
 version of that work for locally conformally flat $C^{2, \ga}$
 manifolds (or at least for the case of doubling of the
 manifold), which will be verified below.
  Then the rest of proof follows from that in \cite{LL03} as we
 explain later.

 Let $(M^n, g)$ be a locally conformally flat compact manifolds
 with totally geodesic boundary. We denote a boundary neighborhood
 in $M$ by $U_a \cup \de' U_a$ where $U_a$ is open and $\de' U_a = \de
 M \cap \de U_a$ is a segment on the boundary. By definition,
 there is a conformal map $\phi_a: U_a \cup \de' U_{a} \rightarrow
 V_a \cup \de'V_a \subset S^n_+ \cup S^{n-1}$ such that $V_a \subset S^n_+$
 and $\de' V_a$ is on the equator. Denote the doubling of $M$ by
 $N = M \cup M^*.$ We will define a locally conformally flat smooth
 structure on $N$. Define the corresponding conformal map $\phi_a^*$
 from $U_a^* \subset M^*$ to $V_a^* \subset S^n_-$ through reflection.
 If $\phi_b$ and $\phi_b^*$ is another pair of conformal map such
 that $U_a \cap U_b$ (and thus $U_a^* \cap U_b^*$) is nonempty,
 then there is a conformal transformation $\Phi$ from $\phi_a (U_a \cap
 U_b)$ to $\phi_b (U_a \cap U_b).$ Similarly, there is  a
 corresponding conformal transformation $\Phi^*$ on the counterpart.
 By Liouville theorem, the conformal transformations $\Phi$ and $\Phi^*$ can
 be extended to  conformal transformations on $S^n,$ still denoted by
 $\Phi$ and $\Phi^*.$  If we can prove that $\Phi = \Phi^*,$
 then they define a locally conformally flat smooth structure on
 $N.$ Suppose that $\Phi$ and $\Phi^*$ are not equal. Then $\Phi^{-1} \circ
 \Phi^*$ is not the identity map on $S^n.$ Notice that it is the
 identity map on $\phi_a (\de' (U_a \cap U_b)),$ which is a
 co-dimensional one submanifold contained in the equator. Thus,
  $\Phi^{-1} \circ \Phi^*$ must be a reflection with respect to the equator (see for
 example, Chapter A in \cite{BP92}). This gives us an contradiction because
 $\Phi^{-1}$ can not map $\phi_b^*  (U_a^* \cap U_b^*) \subset
 S^n_-$ to $\phi_a  (U_a \cap U_b) \subset
 S^n_+.$

 We still denote the metric extended to $N$ by $g.$  $(N, g)$
 is then a locally conformally flat closed manifold with $g \in C^{2, \ga}.$
 We also have $(N, g) \ncong (S^n, g_c)$ because $(M, g) \ncong (S^n_+, g_c).$
 Moreover, each side of differentiations in $g$ is defined. We
 can  follow the proof in \cite{SY88} to show that there is a $C^{2, \ga}$
 developing map from the universal cover $\widetilde{N}$ to $S^n.$
 Note that each side of third derivatives in $g$ is defined.
 Hence, the Liouville theorem is still valid since the proof is by an
 ordinary-differential-equations
 approach. Now that $R_g > 0$ on $\widetilde{N},$
 by the same argument as in \cite{SY88}, the developing map is
 injective.  Solutions on $M$ to (\ref{e:main'}) are extended naturally to
 the ones in $C^{2, \ga}$ on $\widetilde{N}.$ To get $C^0$ bounds of $u,$ the
  proof follows from that in \cite{LL03} (proof of (1.44))
 with some revise as we state below. First, instead of using Theorem 1.20 in \cite{LL03},
  we use local estimates in \cite{Chen05} to drop the condition $H_1$ in
 establishing (4.1) in \cite{LL03}. We also drop condition (1.41) in \cite{LL03} by
 noting that the function $F$ we consider is homogeneous, symmetric
 and normalized with $F(e) = 1.$ After getting lower bounds of $u$
 on $(M, g)$ (or equivalently upper bounds in \cite{LL03} because the functions are chosen
 differently), by local estimates \cite{Chen05} and
 Theorem~\ref{t:bdy} we obtain the Harnack inequality
  $$\max_M u \leq C \min_M u.$$ Thus, we only need to prove that
  $\min_M u$ is upper bounded. This follows from the fact that at
  the minimum point $x_0$, we have $\hat A = \hess u + A_g \geq A_g,$ where we
  use the boundary condition $u_n =0$ when $x_0$ is on the boundary.
  Therefore,
   $$e^{-2 \min _M u} = F( \hess u (x_0) + A_g (x_0)) \geq F( A_g(x_0)) > 0. $$

 \noindent (2) $C^2$ estimates.

  Interior $C^2$ estimates are proved in \cite{Chen05}.
  To get boundary $C^2$ estimates, we use Fermi coordinates in a
  tubular neighborhood $\de M \times [0, \iota]$ of the boundary.
  Note that $\de M$ is compact so $\iota$ is a positive
  number.
   Since $g$ is locally  conformally flat,  in a local chart
  we can choose a flat metric $g_0,$ which is conformal to $g,$ such that
  $\mu_{g_0}$ is a constant and $\hat \mu$ is zero.
  Thus, by Theorem~\ref{t:bdy}, we obtain boundary $C^2$ estimates
  in each half ball $\overline {B}^+_r.$ Since $\de M$ is compact,
  there are finitely many local charts of a tubular neighborhood
  of the boundary. We then get the desired estimates.

 \noindent(3) $C^{\infty}$ estimates.

 Once we have $C^2$ bounds,  $F$ is uniformly elliptic and
 concave. By Evans-Krylov \cite{Evan82} and Lions-Trudinger
 \cite{LT86}, we have $C^{2, \ga}$ estimates in the interior and
 on the boundary, respectively. Higher order regularity follows by
 standard elliptic theory.

\epf

\section{Proof of Theorem~\ref{t:bdy}}\label{s:bdy}
 In this Section, we prove boundary estimates.
 We will use an idea in \cite{Chen05} to
derive boundary $C^2$ estimates directly from boundary $C^0$
estimates.

\bpf
 Since $g$ is flat, by Codazzi equation $\mu_g = \mu $ is constant on $\Sigma_r$.
 Let $\hat A = \hess u + du \otimes du - \frac{1}{2}|\gra u|^2 g.$
 The condition $\Gamma^+_1 \subset \Gamma$
 gives $$0 < \sigma_1(\hat A) = trace\, \hat{A} = \Delta u -\frac{n-2}{2} |\gra u|^2.$$
 Thus, $\Delta u $ is positive and
 \beq \label{i:gra}
 |\gra u|^2 < C \Delta u .
 \eeq\par
 (1) We show that $u_{nnn}$ can be controlled on the boundary.
 Differentiating the equation on both sides in the normal
 direction at a boundary point, we get
 $$ (f e^{-2u})_n = \sum_{\ga, \gb} F^{\ga \gb} \hat{A}_{\ga \gb,n} +
 F^{nn} \hat{A}_{nn,n},$$ where we have used $F^{\alpha n} = 0$ by
 Lemma~\ref{l:bdy}.\\
 For case (a), by Lemma~\ref{l:bdy} again $\hat{A}_{\ga \gb, n} = 2 \mu \hat{A}_{\ga \gb}$. Thus,
  \begin{align}
   (f e^{-2u})_n &= \sum_{\ga, \gb} 2 \mu F^{\ga \gb} \hat{A}_{\ga \gb}
    + F^{nn} \hat{A}_{nn,n} \notag \\
    &= 2 \mu F + F^{nn} (\hat{A}_{nn,n} - 2 \mu \hat{A}_{nn})
   = 2 \mu f e^{-2u} +  F^{nn} (\hat{A}_{nn,n} - 2 \mu \hat{A}_{nn}), \label{i:a}
  \end{align} where the second equality holds by Lemma~\ref{l:sym}
  (a).
 By (\ref{e:nga}) and (\ref{e:ngagb}), we obtain
  \begin{align}
  \hat{A}_{nn,n}- 2 \mu \hat{A}_{nn} &= u_{nnn} +  (u_n- 2 \mu) u_{nn} - \sum_{\ga} u_{\ga} u_{\ga n}
     - 2 \mu ( u^2_n - \frac{1}{2} |\gra u|^2 )\notag \\
      &= u_{nnn} - 3 \mu u_{nn} - \mu^3.\notag
  \end{align}
 Returning to (\ref{i:a}), we get
 $$  -C e^{-2 u} \leq  F^{nn} ( \hat{A}_{nn,n}- 2 \mu \hat{A}_{nn})
     \leq F^{nn} (u_{nnn} - 3 \mu u_{nn} + C).
 $$
  On the other hand, by condition (S3) we have
 $ F^{nn} \geq \epsilon \frac{F}{\gs_1} \geq  \frac{C}{\Delta u} e^{-2u}.$
 Hence, there is a positive number $L$ such that
 \beq \label{i:third}
 u_{nnn} \geq -L \Delta u + 3 \mu u_{nn} - C
 \eeq is true for every point on the boundary, where $L$ and $C$
 depends on $n, \epsilon, \mu, \|f\|_{C^1}$ and $\inf f.$

 For case (b), by Lemma~\ref{l:bdy} (b) we get
  \begin{align}
   (f e^{-2u})_n &= \sum_{\ga, \gb} F^{\ga \gb} ( 2\mu \hat{A}_{\ga
   \gb}- \hat \mu e^{-u} (\hat A_{\ga \gb} + \hat A_{nn} g_{\ga
   \gb}))  + F^{nn} \hat{A}_{nn,n} \notag \\
    &= (2 \mu - \hat \mu e^{-u}) f e^{-2u} - \hat \mu e^{-u} \sum_{\ga}F^{\ga \ga} \hat A_{nn}
     + F^{nn} (\hat{A}_{nn,n} - (2 \mu - \hat \mu
     e^{-u})\hat{A}_{nn}),\notag
  \end{align} where the second equality holds by Lemma~\ref{l:sym}
  (a). Note that $\hat \mu$ is positive. Thus, if $\hat A_{nn} \geq 0,$
  then
   $$ -C e^{-2u}\leq F^{nn} (\hat{A}_{nn,n} - (2 \mu - \hat \mu
   e^{-u})\hat{A}_{nn}).$$
  On the other hand, if $\hat A_{nn} < 0,$ by condition (A) we
  have $$ -C e^{-2u}\leq F^{nn} (\hat{A}_{nn,n} - (2 \mu + \rho \, \hat
  \mu e^{-u})\hat{A}_{nn}),$$ where we drop the term $ F^{nn} \hat \mu
     e^{-u} \hat{A}_{nn}$ since it is negative.
  Hence, in both cases we obtain
   \beq \label{i:a'}
   -C e^{-2u}\leq F^{nn} (\hat{A}_{nn,n} - 2 \mu \hat{A}_{nn} +
   C |\hat{A}_{nn}|).\eeq
  Now  by (\ref{e:nga}) and (\ref{e:ngagb}) and combined with  a basic fact
  that if $\Gamma^+_2 \subset \Gamma,$ then $|u_{ij}| \leq C \Delta u,$
   we get
 $$  \hat{A}_{nn,n}- 2 \mu \hat{A}_{nn} + C |\hat{A}_{nn}|
     \leq u_{nnn} +(- 3 \mu + \hat \mu e^{-u}) u_{nn} + C \Delta u + C.
 $$
 Returning to (\ref{i:a'}), note that by condition (S3) we have
 $ F^{nn} \geq \epsilon \frac{F}{\gs_1} \geq  \frac{C}{\Delta u} e^{-2u}.$
 Hence, there is a positive number $L$ such that
 \beq \label{i:third'}
 u_{nnn} \geq -L \Delta u + (3 \mu - \hat \mu e^{-u}) u_{nn} - C
 \eeq is true for every point on the boundary, where $L$ and $C$
 depends on $n, \epsilon, \rho, \mu, \hat \mu, \inf u, \|f\|_{C^1}$ and $\inf f.$
\vskip 1em
 (2) We will show that $\Delta u$ is bounded. The follow proof is for
 both cases (a) and (b), while the number $C$ is understood as a
 constant depending on $ n, r, \epsilon, \mu, \| f\|_{C^2}$ and $\inf
  f$ for case (a), and $ n, r, \epsilon, \rho, \mu, \hat \mu, \inf u, \| f\|_{C^2}$ and
$\inf  f$ for case (b), respectively.

 Let $H = \eta ( \Delta u +
|\gra u|^2 + n\mu \, u_n) e^{a\,x_n}$ where $a$ is some number
decided later. Denote $r^2 :\equiv \sum_i x^2_i.$ Let $\eta(r)$ be
a cutoff function such that $0 \leq \eta \leq 1$, $\eta = 1$ in
$\overline{B}^+_{\frac{r}{2}}$ and $\eta = 0$ outside
$\overline{B}^+_r,$ and also $|\gra \eta|< C
\frac{{\eta}^{\frac{1}{2}}}{r}$ and $|\hess \eta|< \frac{C}{r^2}.$
By (\ref{i:gra}), $\Delta u$ is positive. Without loss of
generality, we may assume $r= 1$ and
 $$K= \Delta u + |\gra u|^2 + n\mu \, u_n \gg 1.$$
 At a boundary point, note that $\eta_n = 0$ because $\eta = \eta (r).$
 Differentiating $H$ in the normal direction
 and using (\ref{e:nga}) and (\ref{e:ngagb}) gives
 \begin{eqnarray*}
 H_n&=& \eta_n ( K e^{ax_n})+ \eta (K_n  + a K) e^{a x_n}=\eta (K_n + a K) e^{a x_n}  \\
 &\geq& \eta ( ( u_{nnn} + (2 \mu - \hat \mu e^{-u}) K + (- 3 \mu + \hat \mu e^{-u})u_{nn} -C ) + a
 K)e^{a x_n}.
\end{eqnarray*}
By (\ref{i:gra}) and the inequalities (\ref{i:third}) and
(\ref{i:third'}) for
 cases (a) and (b), respectively, we obtain
  \begin{eqnarray*}
 H_n &\geq& \eta ( - L \Delta u + (2 \mu - \hat \mu e^{-u}) K  -C  + a K)e^{a
 x_n}> 0
\end{eqnarray*}
for $a > L - 2\mu + \hat \mu \sup e^{-u} +1.$ Thus, $H$ increases
toward the interior and the maximum of $H$ must happen at some
point $x_0$ in the interior.

At the maximal point $x_0$, we have
\begin{equation}\label{e:star'}
H_i = \eta_i (K e^{a x_n}) + \eta e^{a x_n}(K_i + a K \gd_{in})=0,
\end{equation}and
 \begin{eqnarray*}
 H_{ij} &=& \eta_{ij} ( K+ e^{a x_n}) + \eta_i ( K e^{a x_n})_j +
  \eta_j ( K e^{a x_n})_i + \eta  (K e^{a x_n})_{ij}\\
    &=& (\eta_{ij}- 2 \eta^{-1} \eta_i \eta_j) K e^{a x_n} +
    \eta e^{a x_n}(K_{ij} + a K_i \gd_{jn} + a K_j \gd_{in} + a^2 K
    \gd_{in} \gd_{jn})
  \end{eqnarray*}
is negative semi-definite, where in the second equality we have
used (\ref{e:star'}). Using the positivity of $F^{ij}$ and
(\ref{e:star'}) again to replace $K_i$ and $K_j$, we get
 \begin{align}
 0 \geq F^{ij} H_{ij} e^{-a x_n} &= F^{ij} ((\eta_{ij}- 2 \eta^{-1} \eta_i \eta_j) K +
    \eta (K_{ij} - a \frac{\eta_i}{\eta} K \gd_{jn} - a \frac{\eta_j}{\eta} K
    \gd_{in}- a^2 K \gd_{in} \gd_{jn}))\notag \\
      &\geq   \eta F^{ij} K_{ij} -C \sum_i F^{ii} K, \label{i:b}
 \end{align}
 where we use conditions on $\eta$ in the inequality.
 By direct computations,
 $$
  F^{ij} K_{ij} = F^{ij}u_{llij} + F^{ij}( 2 u_{li}u_{lj} + 2 u_l u_{lij}
  + n \mu u_{nij}) = I + II.$$

 For I, notice that
 $$\hat{A}_{ij,ll} = u_{ijll} + 2 u_{il}u_{jl} + u_i u_{jll} + u_j u_{ill} -
   (u_k u_{kll} + u^2_{kl})g_{ij}.$$
 Then
 $$
  I = F^{ij} (\hat{A}_{ij,ll}- 2 u_{li}u_{lj}- 2 u_{ill} u_j +
  (u_{lk} ^2 + u_k u_{kll}) g_{ij}),$$
 where $F^{ij}(u_i u_{jll})= F^{ij}(u_j u_{ill})$ because $F^{ij}$ is symmetric.
 Now using (\ref{e:star'}) to replace $u_{lli}$ and $u_{kll}$
 yields
  \begin{eqnarray*}
   I &=& F^{ij}\hat{A}_{ij,ll} + F^{ij}( - 2 u_{li}u_{lj} -2 u_j (-2 u_l u_{li}
   -n \mu u_{ni} - \frac{\eta_i}{\eta} K - a K \gd_{in})\\
    & &  + (|\hess u|^2 + u_k (-2 u_l u_{lk}- n\mu u_{nk} - \frac{\eta_k}{\eta} K - a K \gd_{kn}))
   g_{ij}).
  \end{eqnarray*}
 By (\ref{i:gra}) and the conditions on $\eta,$ we have
   \begin{eqnarray*}
   I &\geq& F^{ij}\hat{A}_{ij,ll} + F^{ij}( - 2 u_{li}u_{lj} + 4 u_j u_l u_{li}
     +(|\hess u|^2 - 2 u_k u_l u_{lk}) g_{ij})\\
    & &- C \sum_i F^{ii}\eta^{- \frac{1}{2}}( 1 + |\hess u|^{\frac{3}{2}}).
  \end{eqnarray*}

 For II, we use the formula
 $$\hat A_{ij,l} = u_{ijl} + u_i u_{jl}+ u_j u_{il}- u_k u_{kl} g_{ij}$$
 to obtain
  \begin{eqnarray*}
  II &=& F^{ij}( 2 u_{li}u_{lj} + 2 u_l u_{ijl}+ n \mu u_{ijn})
   = F^{ij} (2 u_{li}u_{lj}+ 2 u_l (\hat A_{ij,l} - 2 u_i
   u_{jl} + u_k u_{kl} g_{ij})\\
   & & + n \mu (\hat A_{ij,n}
   -2 u_i u_{jn} + u_k u_{kn} g_{ij}))\\
   &\geq& F^{ij} (2 u_{li}u_{lj}+ 2 u_l \hat A_{ij,l} - 4 u_i
   u_{jl} u_j + 2 u_k u_{kl} u_l g_{ij}+ n \mu \hat
   A_{ij,n}) - C \sum_i F^{ii}  |\hess u|^{\frac{3}{2}}.
  \end{eqnarray*}

 Combining I and II together, we find that

   \begin{eqnarray*}
  F^{ij} K_{ij} &\geq& F^{ij}\hat{A}_{ij,ll} + F^{ij}( - 2 u_{li}u_{lj} + 4 u_j u_l u_{li}
     +(|\hess u|^2 - 2 u_k u_l u_{lk}) g_{ij})\\
    & &+ F^{ij} (2 u_{li}u_{lj}+ 2 u_l \hat A_{ij,l} - 4 u_i
   u_{jl} u_j + 2 u_k u_{kl} u_l g_{ij}+ n \mu \hat
   A_{ij,n})\\
  & &  - C \sum_i F^{ii} \eta^{-\frac{1}{2}}( 1 + |\hess u|^{\frac{3}{2}}).
  \end{eqnarray*}
 Here is the key step of the proof. Three terms from I cancel out
 three terms from II. Thus, after the cancellations we arrive at
    \begin{eqnarray*}
  F^{ij} K_{ij} &\geq& F^{ij}\hat{A}_{ij,ll} + F^{ij}|\hess u|^2  g_{ij}
    + F^{ij} ( 2 u_l \hat A_{ij,l}+ n \mu \hat A_{ij,n})\\
  & & - C \sum_i F^{ii} \eta^{-\frac{1}{2}}( 1 + |\hess u|^{\frac{3}{2}}).
  \end{eqnarray*}
Now returning to (\ref{i:b}),  applying $\eta$ on both sides
produces
 \begin{eqnarray*}
 0 &\geq& \eta^2  F^{ij} K_{ij} -C \sum_i F^{ii} \eta K\\
   &\geq& \eta^2 F^{ij}\hat{A}_{ij,ll} + \eta^2 F^{ij}|\hess u|^2  g_{ij}
    + \eta^2 F^{ij} ( 2 u_l \hat A_{ij,l}+ n \mu \hat A_{ij,n})\\
   & & - C \sum_i F^{ii} ( 1 + \eta^{\frac{3}{2}} |\hess u|^{\frac{3}{2}}).
  \end{eqnarray*}
 By the concavity of $F$ and Lemma~\ref{l:sym}(a), we have
 $F^{ij} \hat A_{ij,ll} \geq (F^{ij} \hat A_{ij})_{ll}= (f
 e^{-2u})_{ll}.$ Hence,
  \begin{eqnarray*}
 0 &\geq&  \eta^2 \sum_i F^{ii}|\hess u|^2 + \eta^2 (f e^{-2u})_{ll}
    + 2 \eta^2 u_l (f e^{-2u})_l + n \mu \eta^2 (f e^{-2u})_n\\
  & &  - C \sum_i F^{ii} ( 1 + \eta^{\frac{3}{2}} |\hess u|^{\frac{3}{2}})\\
   &\geq& \sum_i F^{ii} ( \eta^2 |\hess u|^2   - C - C \eta |\hess u|- C \eta^{\frac{3}{2}} |\hess
   u|^{\frac{3}{2}}).
  \end{eqnarray*}
 This gives $(\eta |\hess u|)(x_0) \leq C.$ Hence, for $x \in
\overline{B}^+_{\frac{r}{2}},$ we have that $H = ( \Delta u +
|\gra u|^2 + n \mu \, u_n) e^{a\,x_n}$ is bounded. Thus, $\Delta
u$ is bounded. By (\ref{i:gra}), $|\gra u|$ is also bounded.
\vskip 1em
 (3) To get the Hessian bounds, for case (b) it
 follows immediately by the fact that if $\Gamma^+_2 \subset \Gamma,$
 then $|u_{ij}| \leq C \Delta u.$
 As for case (a), note that from (2) above we have $\eta \Delta u <
 C$ and $\eta |\gra u|^2 < C.$ Consider the maximum of $\eta (\hess u + du \otimes du
 + \mu u_n g) e^{a x_n}$ over the set $(x, \xi) \in (B_1^+, S^n).$
 We will show that at the maximum,  $x$ can not belong to the
 boundary. If $\xi$ is in the tangential direction, without loss
 of generality, we can assume $\xi$ is in $e_1$ direction. We have
 \begin{eqnarray*}
 (\eta (u_{11} + u_1^2 + \mu u_n) e^{a x_n})_n &=& \eta (u_{11n} + 2 u_1 u_{1n} +
  \mu u_{nn} + a(u_{11}+ u_1^2+ \mu u_n)) e^{a x_n}\\
    &=& \eta e^{a x_n} ( (2 \mu + a) (u_{11} + u_1^2 + \mu u_n) +
    \mu^3)> 0
 \end{eqnarray*} for $a > -2\mu + 1.$ If $\xi$ is in the normal
 direction, we first have that $\Delta u \leq n(u_{nn} + \mu^2) \leq n u_{nn} + C.$
 By (\ref{i:third}), we obtain
  \begin{eqnarray*}
  (\eta (u_{nn} + u_n^2 + \mu u_n) e^{a x_n})_n &=& \eta (u_{nnn} - \mu u_{nn} +
   au_{nn}) e^{a x_n}\\
  &\geq& \eta e^{a x_n} (-L \Delta u + 2 \mu u_{nn} + a u_{nn}- C)\\
  &\geq&  \eta e^{a x_n} (-n L u_{nn} + 2 \mu u_{nn} + a u_{nn}-
  C)> 0
  \end{eqnarray*} for $a> nL- 2\mu + 1.$ Thus, we  conclude that at the
  maximum, $x$ must be in the interior. We then perform
  similar computations as before using the inequality $\eta |\gra u|^2 < C$
  to get the Hessian bounds. We omit the details here.
 \epf

 Department of Mathematics, Princeton University, Princeton, NJ \par
 Email address: \textsf{szuchen@math.princeton.edu}
 \vskip 1em
 \noindent This article appeared in the American Journal of Mathematics,
 Volume 00, Issue 00, Year, pages 000-000, Copyright c The Johns Hopkins University Press.


\begin{thebibliography}{10}

\bibitem{BP92}
Riccardo Benedetti and Carlo Petronio.
\newblock {\em Lectures on hyperbolic geometry}.
\newblock Universitext. Springer-Verlag, 1992.

\bibitem{CNS-1}
L.~Caffarelli, L.~Nirenberg, and J.~Spruck.
\newblock The {D}irichlet problem for nonlinear second-order elliptic
  equations. {I}. {M}onge-{A}mp\`ere equation.
\newblock {\em Comm. Pure Appl. Math.}, 37(3):369--402, 1984.

\bibitem{CNS-3}
L.~Caffarelli, L.~Nirenberg, and J.~Spruck.
\newblock The {D}irichlet problem for nonlinear second-order elliptic
  equations. {III}. {F}unctions of the eigenvalues of the {H}essian.
\newblock {\em Acta Math.}, 155(3-4):261--301, 1985.

\bibitem{CGY02}
Sun-Yung~A. Chang, Matthew~J. Gursky, and Paul Yang.
\newblock An a priori estimate for a fully nonlinear equation on
  four-manifolds.
\newblock {\em J. Anal. Math.}, 87:151--186, 2002.

\bibitem{CGY02a}
Sun-Yung~A. Chang, Matthew~J. Gursky, and Paul~C. Yang.
\newblock An equation of {M}onge-{A}mp\`ere type in conformal geometry, and
  four-manifolds of positive {R}icci curvature.
\newblock {\em Ann. of Math. (2)}, 155(3):709--787, 2002.

\bibitem{Chen05b}
Szu-yu~Sophie Chen.
\newblock Conformal deformation on manifolds with boundary.
\newblock preprint.

\bibitem{Chen05}
Szu-yu~Sophie Chen.
\newblock Local estimates for some fully nonlinear elliptic equations.
\newblock {\em Int. Math. Res. Not.}, (55):3403-3425, 2005

\bibitem{Es92}
Jos{\'e}~F. Escobar.
\newblock The {Y}amabe problem on manifolds with boundary.
\newblock {\em J. Differential Geom.}, 35(1):21--84, 1992.

\bibitem{Evan82}
Lawrence~C. Evans.
\newblock Classical solutions of fully nonlinear, convex, second-order elliptic
  equations.
\newblock {\em Comm. Pure Appl. Math.}, 35(3):333--363, 1982.

\bibitem{Gar59}
Lars G{\.a}rding.
\newblock An inequality for hyperbolic polynomials.
\newblock {\em J. Math. Mech.}, 8:957--965, 1959.

\bibitem{Gb05}
Bo~Guan.
\newblock Conformal metrics with prescribed curvature functions on manifolds
  with boundary.
\newblock preprint.


\bibitem{GW03}
Pengfei Guan and Guofang Wang.
\newblock Local estimates for a class of fully nonlinear equations arising from
  conformal geometry.
\newblock {\em Int. Math. Res. Not.}, (26):1413--1432, 2003.

\bibitem{GW03a}
Pengfei Guan and Guofang Wang.
\newblock A fully nonlinear conformal flow on locally conformally flat
              manifolds.
\newblock {\em J. Reine Angew. Math.}, 557:219--238, 2003.


\bibitem{GV05}
Matthew~J. Gursky and Jeff~A. Viaclovsky.
\newblock Prescribing symmetric functions of eigenvalues of {S}chouten tensor.
\newblock to apear in {A}nn. of {M}ath.

\bibitem{GV03}
Matthew~J. Gursky and Jeff~A. Viaclovsky.
\newblock A fully nonlinear equation on four-manifolds with positive scalar
  curvature.
\newblock {\em J. Differential Geom.}, 63(1):131--154, 2003.

\bibitem{LL03}
Aobing Li and Yanyan Li.
\newblock On some conformally invariant fully nonlinear equations.
\newblock {\em Comm. Pure Appl. Math.}, 56(10):1416--1464, 2003.

\bibitem{Li89}
Yan~Yan Li.
\newblock Degree theory for second order nonlinear elliptic operators and its
  applications.
\newblock {\em Comm. Partial Differential Equations}, 14(11):1541--1578, 1989.

\bibitem{LT86}
P.-L. Lions and N.~S. Trudinger.
\newblock Linear oblique derivative problems for the uniformly elliptic
  {H}amilton-{J}acobi-{B}ellman equation.
\newblock {\em Math. Z.}, 191(1):1--15, 1986.

\bibitem{Mit70}
D.~S. Mitrinovi{\'c}.
\newblock {\em Analytic inequalities}.
\newblock Springer-Verlag, New York, 1970.

\bibitem{Reilly}
Robert~C. Reilly.
\newblock On the {H}essian of a function and the curvatures of its graph.
\newblock {\em Michigan Math. J.}, 20:373--383, 1973.

\bibitem{SY88}
R.~Schoen and S.-T. Yau.
\newblock Conformally flat manifolds, {K}leinian groups and scalar curvature.
\newblock {\em Invent. Math.}, 92(1):47--71, 1988.

\bibitem{Sch91}
Richard~M. Schoen.
\newblock On the number of constant scalar curvature metrics in a conformal
  class.
\newblock In {\em Differential geometry}, volume~52, pages 311--320. Longman
  Sci. Tech., 1991.


\bibitem{Tr90}
Neil~S. Trudinger.
\newblock The {D}irichlet problem for the prescribed curvature equations.
\newblock {\em Arch. Rational Mech. Anal.}, 111(2):153--179, 1990.

\end{thebibliography}
\end{document}